\documentclass[12pt]{article}
\usepackage{amsmath,amsfonts}
\usepackage{graphicx}
\begin{document}

\baselineskip=24pt

\title{Derivation of a Stochastic Neutron Transport Equation}

\author{Edward J. Allen \\
Department of Mathematics and Statistics \\ Texas Tech University
 \\ Lubbock, Texas 79409-1042 \\
Email: edward.allen@ttu.edu\\}

\date{}

\maketitle

\noindent {\large This paper is to be published in the
 {\it Journal of Difference Equations and Applications.}}

\maketitle

\begin{abstract}

Stochastic difference equations and a stochastic partial
differential equation (SPDE) are simultaneously derived
 for the time-dependent neutron angular density in a general
 three-dimensional medium where the neutron angular density is a
 function of position, direction, energy, and time.
 Special cases of the equations are given such as transport in one-dimensional
 plane geometry with isotropic scattering and transport
 in a homogeneous medium. The stochastic equations are derived from basic
 principles, i.e., from the changes that occur in a small time interval.
 Stochastic difference equations of the neutron angular density
 are constructed, taking into account the
 inherent randomness in scatters, absorptions, and source neutrons.
 As the time interval decreases, the stochastic difference
 equations
 lead to a system of It\^{o} stochastic differential equations (SDEs). As
 the energy, direction, and position intervals decrease, an SPDE
 is derived for the neutron angular density.  Comparisons
 between numerical solutions of the stochastic difference
 equations and independently formulated Monte Carlo calculations
 support the accuracy of the derivations.\bigskip

 \end{abstract}

\noindent {\bf Short running title:} Stochastic Neutron Transport



\newpage

\section{Introduction}
\label{intro}

  In the present investigation,
  stochastic versions of the deterministic neutron transport equation
   are derived. Specifically, stochastic difference equations and a stochastic
    partial differential equation (SPDE) are simultaneously
    derived that account for the random
    effects of absorptions, scatters, and source particles and
    generalize the standard deterministic neutron transport
    equation. Numerical approximations of the SPDE, through solution of the system
    of stochastic difference equations, provide
    approximations to the randomly varying neutron densities and yield insight
    into the random behavior of neutron transport.
    The stochastic transport equations are
 most useful for problems involving low numbers of neutrons. As
 the coefficient of variation is often approximately inversely
 proportional to the square root of the population size, the deterministic and stochastic
 transport equations yield essentially the same results for high
 numbers of neutrons.

        There are alternate but apparently equivalent ways to
 derive a system of stochastic differential equations (SDEs) for a
 randomly varying dynamical problem. The first way involves deriving a master
 equation for the random process \cite{gillespie, vankampen}.
 A master equation
 is a differential form of the Chapman-Kolmogorov equation
 involving transition probabilities and is a probability conservation
 equation for the probabilities of separate states. If the transition densities
 are expanded in a parameter that defines the size of the fluctuations
 or jumps, then a Fokker-Planck equation or forward Kolmogorov equation is
 obtained in the first few terms of the expansion. As the
 probability density of an SDE system satisfies a certain forward
 Kolmogorov equation, this procedure infers a particular SDE system that
 approximates the random dynamics of the problem. In this
 procedure, a system of stochastic difference equations is not
 derived as an intermediate step. A second way to
 derive a system of SDEs for a randomly varying problem is by
 studying the changes in the process for a short time interval
 which gives a discrete stochastic model. The discrete stochastic model
 infers a system of stochastic difference equations which, in turn, leads
 to an appropriate SDE system.
 For example, consider a randomly
 varying problem
 where $\vec X =[X_1, X_2, \dots, X_N]^T$
 is a random vector of $N$ components for the
 problem.  Let $\Delta \vec X$ be the change in the process for a
 small time interval $\Delta t$. The expectations
 $\vec \mu(\vec X,t) = E(\Delta \vec X)/\Delta t $ and
 $V(\vec X, t) = E(\Delta \vec X (\Delta \vec X)^T)/\Delta t$
 are determined and a stochastic difference equation approximation
 for the problem has
 the form:
 \begin{equation}
\vec X(t+\Delta t) =  \vec X(t)+ \vec \mu(\vec X(t),t) \Delta t +
(V(\vec X(t),t))^{1/2} \; \sqrt{\Delta t} \; \vec \eta_t
\label{eq1.4a}
\end{equation}
where $\vec \eta_t$ is a vector of length $N$ of independent
normally distributed random numbers with zero mean and unit
variance. Finally, the SDE system that approximates the behavior
of the
 randomly varying process has the form
\begin{equation}
d \vec X(t) =  \vec \mu(\vec X,t) dt + (V(\vec X,t))^{1/2} \; d
\vec W(t) \label{eq1.4b}
\end{equation}
where $\vec W(t)$ is a vector of length $N$ of independent Wiener
processes. It can be shown that the probability density of
solutions of the stochastic system (\ref{eq1.4a}) or
(\ref{eq1.4b}) approximates the probability density of the
original randomly varying process
\cite{eallenbook,allenallenetal,allen03}. Assume now that there
are $M$ possible changes in the process with probabilities $p_j
\Delta t$ for $j=1,2, \dots, M$ for small $\Delta t$.  Also,
assume that the $j$th change alters the $i$th component by
$\lambda_{j,i}$. Then, the elements of $\vec \mu$ are given by
$\mu_i = \sum_{j=1}^{M} p_j \lambda_{j,i}$  and the elements of $N
\times N$ matrix $V$ are given by $v_{i,l}= \sum_{j=1}^{M} p_j
\lambda_{j,i} \lambda_{j,l}$ for $i,l=1,2, \dots N$. In addition,
an equivalent SDE system to (\ref{eq1.4b}) is
\begin{equation}
d \vec X(t) =  \vec \mu(\vec X,t) dt + C(\vec X,t) \; d \vec
W^*(t) \label{eq1.5}
\end{equation}
where $N \times M$ matrix $C$ has elements $c_{i,j} =
\lambda_{j,i} p_j^{1/2}$ and $\vec W^*(t)$ is a vector of length
$M$ of independent Wiener processes. Furthermore, equation
(\ref{eq1.5}) is also obtained from the master-equation approach
under the same assumptions \cite{gillespie}. Thus, the two
derivation procedures produce identical SDE systems for this
general $N$-component and $M$-change process.

   These two derivation procedures produce very reasonable stochastic
 equation models for a given phenomenon. For randomly varying systems
 where the dependent variable depends
on time and on secondary independent variables, a stochastic
partial differential equation (SPDE) may be derived by
 replacing the Wiener processes in the SDE system with
appropriate Brownian sheets and letting the intervals in the
remaining independent variables go to zero. The resulting equation
is an SPDE model for the phenomenon.

     In this paper, stochastic difference equations and a stochastic partial differential equation
are derived for the transport of neutrons in matter. In neutron
transport,
   captures, scatters, fissions, and source neutrons occur randomly.
   As a result, the neutron angular density varies stochastically.
   The relative magnitude of the random behavior of the angular density
   is pronounced for low neutron densities, such as during reactor startup,
   but decreases as the neutron density increases.
    The standard or deterministic neutron transport equation
(or Boltzmann neutron transport equation) describes the expected
or probable neutron angular density with respect to position,
direction, energy, and time  \cite{bellglasstone}. Solutions of
the deterministic neutron transport equation
 provide average values of the neutron angular density; actual realizations
 with time of the neutron angular densities, that include random effects
   from neutron interactions and sources, are not obtained.

      A stochastic partial
differential equation is derived in this paper for neutron
transport in a general three-dimensional
      absorbing and anisotropic-scattering medium where the neutron angular density
      depends on position, direction, energy, and time.
 In the present investigation, the medium is assumed to be constant
 with respect to material composition, i.e., zero power noise.
    Special random effects, for example, from randomly varying boundary conditions
    or from a medium that is randomly varying \cite{williamslarsen,williams}
    are not studied in the present investigation although
    generalizations of the SPDE to approximate such conditions
    may be possible. Using the derived stochastic
   neutron transport equation,
 sample paths (realizations) of the randomly varying neutron angular
 densities can be approximately computed. In addition, after computing
 many sample
 paths, moments of the neutron densities, for example, can be estimated.
    The stochastic neutron transport equation is derived from basic
   principles, i.e., from the changes in the system that occur in a
   small time interval. The dynamical system is
  studied to determine the different independent random
  changes that occur. Appropriate terms are
  identified for these changes in developing a  stochastic
  difference system where all independent variables are discrete. As the time interval
  goes to zero, a certain stochastic differential equation
  (SDE) system is inferred (e.g.,
  \cite{eallenbook,allenallenetal,allen03,hayesallen,sharpallen}).
  Next, multidimensional Brownian sheets replace
  the Wiener processes. As the intervals in the remaining
  independent variables go to zero, the SDE system leads to an
  SPDE  \cite{eallen0708,eallenjbd}.
  It is illustrated how the stochastic transport equation can be solved
  computationally through numerical solution of a stochastic
  difference system.

    The neutron transport equation is of fundamental importance in
    nuclear reactor theory and shielding design \cite{bellglasstone,duderstadtmartin,hetrick,lewismiller}.
     The stochastic nature of the neutron
     transport process has been of interest for many years.
     Classic studies of the stochastic theory of neutron transport
     are given in \cite{bella,bellb}. In particular,
     let $p_n(R,t_f; \vec x,\vec v,t)$
     be the probability that a neutron with position $\vec x$ and velocity
     $\vec v$ at time $t$ leads to $n$ neutrons in region $R$ of
     $\vec x, \vec v$ space at time $t_f$. In \cite{bellb},
     a non-linear integro-differential equation for the
     probability generating function is derived for
     $p_n(R,t_f; \vec x,\vec v,t)$. The equations
     derived are interesting but complicated and difficult to apply.
     More recently, a master equation approach was used
     to estimate the temporal evolution of the number of neutrons
     in time-varying multiplying systems \cite{palb,pala}.  This approach gives,
     for example, moments of the number of neutrons in the system with
     time. However, a stochastic difference system approximation of neutron transport
     is not determined and, as a result, sample paths of the randomly varying
     neutron densities with respect to energy, position, and
     direction are not estimated.

   In the next section, stochastic difference equations and a stochastic partial differential
   equation are derived for neutron transport in general
   three-dimensional xyz-geometry. The
   changes due to absorptions, fissions, and scatters, which occur
   randomly with probability proportional to the neutron angular density and
   to the time interval, are carefully considered in deriving the
   equations.
    In the third section, several special but useful cases of
   the stochastic neutron transport equation are described such as
   one-dimensional slab geometry with isotropic scattering. In the
   fourth section, stochastic difference equations
   are solved computationally for the randomly varying neutron
   densities and compared with Monte Carlo calculations. The Monte Carlo
   calculational procedures differ considerably from the numerical solution of the
   stochastic transport equations. In the Monte Carlo calculations, the
   dynamical system is checked at each small interval of time to take
   into account scatters, absorptions, and movements for individual
   neutrons. Comparisons between the two different computational
   methods are in close
   agreement indicating that the stochastic difference equations and the
   stochastic partial differential equation
   accurately model the random behavior of neutron transport.
   Although the stochastic equations cannot exactly model the random
   behavior of neutron transport since, for example, the number of neutrons in any
   region is not integer-valued in the stochastic difference system or in the SPDE model,
   the calculations
   indicate that the stochastic neutron transport model is accurate. In addition, the
   stochastic transport equations provide insight into the random dynamics of neutron
   transport and can be efficiently solved computationally using
   the stochastic difference equations.
   Several applications of the stochastic neutron transport equation are
   discussed in the fifth section before the investigation is summarized in the
   final section.

\section{Derivation of Stochastic Neutron Transport Equations}
\label{der}

The neutron transport equation in xyz-geometry can be written in
the integro-differential form
\cite{bellglasstone,duderstadtmartin}:
\begin{eqnarray}
&& \dfrac{\partial  \bar N(x,y,z,\mu,\phi,E,t)}{\partial t} = -v
\mu_x \dfrac{\partial \bar N}{\partial x} -v \mu_y
\dfrac{\partial \bar N}{\partial y} -
v \mu_z \dfrac{\partial \bar N}{\partial z} \label{eq2.1} \\
&& + \, Q(x,y,z,E,t) -v \sigma(x,y,z,E) \bar N(x,y,z,\mu,\phi,E)
\nonumber
\\ && + \int_0^{E_{max}} \int_{-1}^1 \int_0^{2 \pi} v'
\sigma' f(x,y,z,\mu',\phi',E', \mu, \phi,E) \bar N' \; d\phi'
d\mu' dE' \nonumber
\end{eqnarray}
 for $(x,y,z) \in ([0,x_{max}] \times [0,y_{max}] \times [0,z_{max}])
  \subset \mathbb R^3$, $-1 \le \mu \le 1$, $0
\le \phi \le 2 \pi$, $0 \le E \le E_{max}$ and $t \ge 0$ where
$\bar N =\bar N(x,y,z,\mu, \phi,E, t)$ is the expected neutron
angular density with respect to position $(x,y,z)$, direction
$\vec \Omega = (\mu, \phi)$, energy $E$, and time $t$ per unit
volume per unit solid angle per unit energy. In (\ref{eq2.1}),
$\sigma' = \sigma(x,y,z,E')$ and $ \bar N' = \bar
N(x,y,z,\mu',\phi',E',t)$. The notation used here is generally
consistent with the notation used in \cite{bellglasstone}. In
particular, $v$ is the neutron speed, $\sigma$ is the total
macroscopic cross section, $\sigma(x,y,z,E')
f(x,y,z,\mu',\phi',E',\mu,\phi,E) \Delta \mu \Delta \phi \Delta E
$ is the probability of a neutron transfer from direction
$(\mu',\phi')$ and energy $E'$ to solid angle $\Delta \mu \Delta
\phi$ about direction $(\mu,\phi)$ with energy $\Delta E$ about
energy $E$. Note that $$\sigma' f(x,y,z,\mu',\phi',E',\mu,\phi,E)
= \sum_{r} \sigma_r' f_r(x,y,z,\mu',\phi',E',\mu,\phi,E) $$ where
$\sigma_r'= \sigma_r(x,y,z,E')$ and the sum includes the separate
interactions $r$ in which neutrons are produced such as elastic
scattering or fission. The parameters $\mu_x$, $\mu_y$ and $\mu_z$
are the direction cosines for the $x$, $y$, and $z$ axes,
respectively. In particular, $\mu_x = \mu = \cos \theta$, $\mu_y =
\sin (\cos^{-1} \mu) \cos \phi$, and $\mu_z = \sin (\cos^{-1} \mu)
\sin \phi$. Thus, $\vec \Omega \cdot \vec \nabla \bar N = \mu_x
\dfrac{\partial \bar N}{\partial x} + \mu_y \dfrac{\partial \bar
N}{\partial y} +\mu_z \dfrac{\partial \bar N}{\partial z}$.
Furthermore, $c(x,y,z,E) = \int_0^{E_{max}} \int_{-1}^1 \int_0^{2
\pi} f(x,y,z,\mu,\phi,E, \mu', \phi',E') \; d\phi' d\mu' dE'$ is
the mean number of neutrons emerging per collision of neutrons of
energy $E$ at position $(x,y,z)$. Finally, $Q(x,y,z,E)$ is the
number of source neutrons per unit solid angle per unit volume per
unit energy and is assumed to be isotropic.

   To simplify the derivation, it is useful to define several other
quantities.  Let $\sigma_c$ be the capture cross section, i.e., the
sum of all the cross sections involving a pure capture event such as
those due to $(n,\gamma), (n,p)$, or $(n,\alpha)$ collisions.  Let
$\hat \sigma(x,y,z,E) = \sigma(x,y,z,E)- \sigma_c(x,y,z,E)$ be the
macroscopic cross section for all interactions other than pure
capture interactions. Let $\hat f$ be defined by the expression
$$ \hat \sigma(x,y,z,E') \hat f(x,y,z,\mu',\phi',E',\mu,\phi,E) =
\sigma' f(x,y,z,\mu',\phi',E',\mu,\phi,E) $$ and define $\hat
c(x,y,z,E)  = \int_0^{E_{max}} \int_{-1}^1 \int_0^{2 \pi} \hat
f(x,y,z,\mu,\phi,E, \mu', \phi',E') \; d\phi' d\mu' dE'$ as the
mean number of neutrons emerging per non-capture collision of
neutrons of energy $E$ at position $(x,y,z)$.

 Equation
(\ref{eq2.1}) is deterministic and random variations in the
neutron angular density due to the inherent randomness in
absorptions, scatters, and fissions cannot be accurately studied
using this equation. To derive a stochastic partial differential
equation generalization of (\ref{eq2.1}), the changes which occur
in the angular density for a small time interval are determined
taking into account interactions and transport.
 A discrete stochastic model of the neutron angular density
 is then constructed which infers a system of stochastic difference
 equations. As the time interval decreases,
 the stochastic difference system leads to a system of It\^{o} stochastic
 differential equations. As the intervals in position, direction, and
 energy decrease, a stochastic partial
 differential equation is derived for the neutron transport
 process.

\subsection{Several Properties of Brownian Sheets}

\begin{figure}[ht]
\begin{center}
\includegraphics[height=2.2in,width=4in]{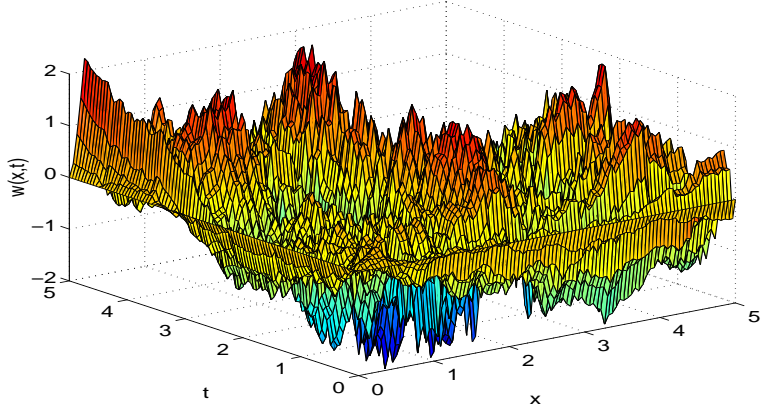}
\caption{A Brownian sheet on $[0,5] \times [0,5]$.}\label{figbs}
\end{center}
\end{figure}

 Before deriving these stochastic equations, it is useful to
 consider
 several properties of Brownian sheets \cite{allennovoselzhang,cabana,walsh}.
 A Brownian sheet on $[0,5] \times [0,5]$ is illustrated in Fig.
\ref{figbs}.
The Brownian sheet $W(x,t)$ satisfies:
\begin{equation*}
\int_t^{t+\Delta t} \int_{x}^{x+\Delta x} \dfrac{\partial^2 W
(x',t')}{\partial t' \partial x'} \, d x' \, d t' \sim {\mathcal
N}(0, \Delta x \Delta t). \label{eq2.1a}
\end{equation*}
That is, the Brownian sheet is independent and normally
distributed over rectangular regions. In addition, if $x_j= j
\Delta x$ for $j=0,1, \dots, J$, where $\Delta x =x_{max}/J$, then
the Brownian sheet defines for $j=1,2, \dots J$, the standard
Wiener processes, $W_j(t)$, where
$$  \dfrac{d W_j(t)}{dt} = \dfrac{1}{\sqrt{\Delta x}}
\int_{x_{j-1}}^{x_{j}}
\dfrac{\partial^2 W (x',t)}{\partial t
\partial x'} \, d x'. $$ Notice that if $t_i=i \Delta t$ for
$i=0,1, \dots,M$, then
$$   \int_{t_{i-1}}^{t_{i}} d W_j(t') = \sqrt{\Delta t}
\, \eta_{i,j}$$ where $\eta_{i,j} \sim {\mathcal N} (0,1)$ for
each $j=1,2, \dots J$ and $i=1,2, \dots M$. Also, standard Wiener
processes can be defined using, for example, three-dimensional
Brownian sheets letting $$  d W_{j,k}(t)= \dfrac{1}{\sqrt{\Delta x
\Delta y}} \int_{x_{j-1}}^{x_{j}} \int_{y_{k-1}}^{y_{k}}
\dfrac{\partial^3 W (x',y',t)}{\partial t \partial y' \partial x'}
\, d y' dx' dt $$ where $W_{j,k}(t)$ is a Wiener process for each
$j$ and $k$. However, notice that $$ \dfrac{\partial^3 W
(x,y,t)}{\partial t
\partial y \partial x} = \lim_{\Delta x, \Delta y \to 0}
\dfrac{1}{\Delta x \Delta y} \int_{x}^{x+\Delta x} \int_{y}^{y+
\Delta y} \dfrac{\partial^3 W (x',y',t)}{\partial t
\partial y' \partial x'} \, d y' dx'. $$

   Finally, it is useful to note that $W(x,t) \ne W_1(x)W_2(t)$
   where $W_1(x)$ and $W_2(t)$ are independent Wiener processes.
 To see this, let $A=[a,b] \times [c,d]$ and $|A| =(d-c)(b-a)$.
 Then,
  $${\displaystyle W(A) = \int_a^b \int_c^d \dfrac{\partial^2 W (x,t)}{\partial t
\partial x} \; dx dy = W(b,d)-W(b,c)-W(a,d)+W(a,c) \sim \mathcal N(0,|A|).}$$
In particular, $W(A)= \sqrt{|A|} \; \eta_{xt}$ where $\eta_{xt}
\sim \mathcal N(0,1)$.  However,
$$ W_{1,2}(A) = \int_a^b \int_c^d d W_1(x) d W_2(t) =
(W_1(b)-W_1(a))(W_2(d)-W_2(c))= \sqrt{|A|} \; \eta_{x}\eta_{t}$$
where $\eta_{x} \sim \mathcal N(0,1)$
 and $\eta_{t} \sim \mathcal N(0,1)$. So, for example, $E((W(A))^4) = 3 |A|^2 $
 whereas $E((W_{1,2}(A))^4) = 9 |A|^2 $.

\subsection{Derivation of stochastic neutron transport equations}

To derive a stochastic neutron transport equation, the changes
which occur in the neutron angular density for a small time
interval $\Delta t$ at time $t_p$ are considered where
$t_p=(p-1)\Delta t$ for $p=1,2, \dots.$ To facilitate finding
these changes, the variables position, direction, and energy are
made discrete. Three-dimensional space is discretized into
rectangular parallelepipeds of length, width, and height $\Delta
x, \Delta y, \Delta z$, respectively. The direction variables
$\mu$ and $\phi$ are discretized with intervals of size $\Delta
\mu =2/L$ and $\Delta \phi =2 \pi/M$, and energy is discretized
into intervals of size $\Delta E =E_{max}/G$. Furthermore, $x_i=
(i-1) \Delta x$ for $i=1,2,\dots, I+1$, $y_j= (j-1) \Delta y$ for
$j=1,2,\dots, J+1$, $z_k= (k-1) \Delta z$ for $k=1,2,\dots, K+1$,
$\mu_l= -1 + (l-1) \Delta  \mu $ for $l=1,2,\dots, L+1$, $\phi_m=
(m-1) \Delta \phi$ for $m=1,2,\dots, M+1$, and $E_g = (g-1) \Delta
E$ for $g=1,2,\dots, G+1$. Now let
$n(x_i,y_j,z_k,\mu_l,\phi_m,E_g,t_p) =
N(x_i,y_j,z_k,\mu_l,\phi_m,E_g,t_p) \Delta x \Delta y \Delta z
 \Delta \mu \Delta \phi \Delta E $ be the number of neutrons
 in the parallelepiped moving in direction $\mu_l,\phi_m$ with energy
 $E_g$. There are several possible changes that can occur to
 $n=n(x_i,y_j,z_k,\mu_l,\phi_m,E_g,t_p)$ in the small time interval $\Delta t$. A
 capture or fission can occur, a neutron can enter or leave one of
 the six faces of the parallelepiped, or a scatter can occur
 resulting in a loss or gain of one neutron. The possible changes
 $\Delta n$ along with their probabilities are listed in Table 1 
 for a small time interval $\Delta t$. Notice that position changes
 occur deterministically, i.e., the neutron position
 is determined by the neutron velocity and, thus, the number of
 neutrons moving from one parallelepiped into an adjacent
 parallelepiped is calculated based on the fraction of neutrons
 crossing the parallelepiped boundary in time $\Delta t$.
In Table 1, the probabilities $p_c, p_{tr1}, p_{tr2}$, and $p_Q$
are given by
\begin{eqnarray*}
 p_c &=& v_g \sigma_c n(x_{i},y_j,z_k,\mu_l,\phi_m,E_g,t_p) \Delta
 t \\
  p_{tr1} &=& v_g \hat \sigma \hat f n \Delta \mu \Delta \phi \Delta E \Delta
 t, \\
  p_{tr2} &=& v_{g'} \sigma' f' n' \Delta \mu \Delta \phi \Delta E \Delta
  t, \;\; \text{and}\\
   p_Q &=&  Q(x_i,y_j,z_k,E_g,t_p) \Delta x \Delta y \Delta z
   \Delta \mu \Delta \phi \Delta E \Delta t.
\end{eqnarray*}
  For example, in Table 1,
$v_g \sigma_c(x_i,y_j,z_k,E_g)
n(x_{i},y_j,z_k,\mu_l,\phi_m,E_g,t_p) \Delta t$ is the probability
that a neutron in the packet undergoes a capture. However, the two
interaction terms, i.e., the term involving transfer into the
packet and the term involving transfer out of the packet, need to
be considered at each position for all the different directions
and energies. Specifically,
$$ p_{tr1} = v_g \hat \sigma
 \hat f(x_{i},y_j,z_k,\mu_l,\phi_m,E_g,\mu_{l'},\phi_{m'},E_{g'})
   n \Delta \mu \Delta \phi \Delta E  \Delta t$$  is the probability
   in time $\Delta t$ that $1/\hat c(x_i,y_j,z_k,E_g)$ neutrons are lost from
   direction
    $\mu_l,\phi_m$ at energy $E_g$ (from the
   packet) when one neutron emerges in direction $\mu_{l'},\phi_{m'}$ with energy
 $E_{g'}$ for each value of $l'$, $m'$, and $g'$
 and $$ p_{tr2} = v_{g'} \sigma'
    f(x_{i},y_j,z_k,\mu_{l'},\phi_{m'},E_{g'},\mu_l,\phi_m,E_g)
   n' \; \Delta \mu \Delta \phi \Delta E \Delta t$$
   is the probability
   that one neutron emerges in direction $\mu_l,\phi_m$ at energy $E_g$ (into the
   packet) when $1/\hat c(x_i,y_j,z_k ,E_{g'})$ neutrons are lost from
    direction $\mu_{l'},\phi_{m'}$ and energy $E_{g'}$ for each value
    of $l'$, $m'$, and $g'$.
 Note, for convenience in Table 1, 
  $n =n(x_i,y_j,z_k,\mu_{l},\phi_{m},E_{g},t_p)$, $n' =
n(x_i,y_j,z_k,\mu_{l'},\phi_{m'},E_{g'},t_p)$, $\sigma=
\sigma(x_i,y_j,z_k,E_g)$, $\sigma' = \sigma(x_i,y_j,z_k,E_{g'})$,
$\sigma_c = \sigma_c(x_i,y_j,z_k,E_{g}),\;$ $\hat \sigma = \hat
\sigma(x_i,y_j,z_k,E_{g}),\;$ $\hat f= \hat
f(x_{i},y_j,z_k,\mu_l,\phi_m,E_g,\mu_{l'},\phi_{m'},E_{g'})$, $f'=
f(x_{i},y_j,z_k,\mu_{l'},\phi_{m'},E_{g'},\mu_l,\phi_m,E_g)$, and
$\hat c= \hat c (x_i,y_j,z_k,E_g)$. Also, it is assumed that the
neutron source, $Q$, is a Poisson process with the probability of
adding one source neutron to the packet in a small time interval
$\Delta t$ equal to $Q(x_i,y_j,z_k,E_g,t_p) \Delta x \Delta y
\Delta z \Delta \mu \Delta \phi \Delta E \Delta t$.

  Table 1 
    defines a discrete stochastic
  model for the neutron transport system. Using these changes and
  probabilities and letting
  the time interval $\Delta t$ approach zero, a system
  of It\^{o} stochastic differential equations can be formulated for
  the dynamics of this random transport process. First, a deterministic
  equation for the expected number of neutrons at time $t_p+ \Delta t$
  can be derived using the
  results of the Table 1.  
   This equation,
  for $\mu_x,\mu_y,\mu_z >0$, is given by:
\begin{eqnarray}
&& \bar n(x_i,y_j,z_k,\mu_l,\phi_m,E_g,t_p+\Delta t)= \bar
n(x_i,y_j,z_k,\mu_l,\phi_m,E_g,t_p)  \label{eq2.2} \\ && - v_g
\sigma_c \bar n(x_{i},y_j,z_k,\mu_l,\phi_m,E_g,t_p) \Delta t
\nonumber
\\ && + \mu_x v_g \Delta t (\bar n(x_{i-1},y_j,z_k,\mu_l,\phi_m,E_g,t_p)
- \bar n(x_{i},y_j,z_k,\mu_l,\phi_m,E_g,t_p))/\Delta x \nonumber \\
&& + \mu_y v_g \Delta t (\bar
n(x_i,y_{j-1},z_k,\mu_l,\phi_m,E_g,t_p)  - \bar
n(x_{i},y_j,z_k,\mu_l,\phi_m,E_g,t_p)) / \Delta y \nonumber
\\ && + \mu_z v_g \Delta t (\bar n(x_i,y_j,z_{k-1},\mu_l,\phi_m,E_g,t_p)
 - \bar n(x_{i},y_j,z_k,\mu_l,\phi_m,E_g,t_p)) /
\Delta z \nonumber \\
&& + Q(x_i,y_j,z_k,E_g,t_p) \Delta x \Delta y \Delta z \Delta \mu
\Delta \phi \Delta E \Delta t
 \nonumber \\
&& - \sum_{l'=1}^L \sum _{m'=1}^M \sum _{g'=1}^G \dfrac{1}{\hat c}
v_g \hat \sigma \hat
f(x_{i},y_j,z_k,\mu_l,\phi_m,E_g,\mu_{l'},\phi_{m'},E_{g'})
   \bar n \Delta \mu \Delta \phi \Delta E \Delta t
    \nonumber \\ && + \sum_{l'=1}^L \sum _{m'=1}^M \sum _{g'=1}^G v_{g'} \sigma'
f(x_{i},y_j,z_k,\mu_{l'},\phi_{m'},E_{g'},\mu_l,\phi_m,E_g)
   \bar n' \Delta \mu \Delta \phi \Delta E \Delta t
   \nonumber
\end{eqnarray}
where $\bar n$ is the expected number of neutrons in the packet.
Letting $\bar n = \bar N \Delta x \Delta y \Delta z \Delta \mu
\Delta \phi \Delta E$ and allowing $\Delta t$ to approach zero as
well as $\Delta x$, $\Delta y$, $\Delta z$, $\Delta \mu$, $\Delta
\phi,$ and $\Delta E$, it is straightforward to show that
(\ref{eq2.2}) yields the standard neutron transport equation
(\ref{eq2.1}). However, the changes and probabilities given in
Table 1 
can be used to derive an SDE model. Indeed, the discrete
stochastic model and the SDE model will have approximately the
same covariance terms as well as mean terms for small $\Delta t$.

The derivation procedure, described in the Introduction for
obtaining equation (\ref{eq1.5}), is now applied using the changes
and probabilities given in Table 1 to obtain the stochastic terms
in the equations. Specifically, for the $i$th component (packet)
of the system, the coefficient of the independent Wiener process
corresponding to the $j$th change is equal to the product of
$\lambda_{j,i}$ with the square root of the probability for the
$j$th change, recalling that $\lambda_{j,i}$ is the amount that
the $j$th change alters the $i$th component. The changes and
probabilities given in Table 1 imply, for $\mu_x,\mu_y,\mu_z > 0$,
that a very reasonable approximation to the discrete stochastic
model satisfies the stochastic difference system:
\begin{eqnarray}
 && n(x_i,y_j,z_k,\mu_l,\phi_m,E_g,t+\Delta t) =
 n(x_i,y_j,z_k,\mu_l,\phi_m,E_g,t) \label{eq2.3a}
 \\ && + \mu_x v_g (n(x_{i-1},y_j,z_k,\mu_l,\phi_m,E_g,t)
  - n(x_{i},y_j,z_k,\mu_l,\phi_m,E_g,t))\Delta t/  \Delta x \nonumber \\
&& + \mu_y v_g (n(x_i,y_{j-1},z_k,\mu_l,\phi_m,E_g,t)
 - n(x_{i},y_j,z_k,\mu_l,\phi_m,E_g,t))\Delta t
/\Delta y \nonumber
\\ && + \mu_z v_g (n(x_i,y_j,z_{k-1},\mu_l,\phi_m,E_g,t)
 -  n(x_{i},y_j,z_k,\mu_l,\phi_m,E_g,t))\Delta t
/\Delta z \nonumber \\ && + Q(x_i,y_j,z_k,E_g,t) \Delta x \Delta y
\Delta z \Delta \mu \Delta \phi \Delta E \Delta t - v_g \sigma_c
n(x_{i},y_j,z_k,\mu_l,\phi_m,E_g,t) \Delta t \nonumber \\
 && - \sum_{l'=1}^L
\sum _{m'=1}^M \sum _{g'=1}^G \dfrac{1}{\hat c} v_g \hat \sigma
\hat f(x_{i},y_j,z_k,\mu_l,\phi_m,E_g,\mu_{l'},\phi_{m'},E_{g'})
   n \Delta \mu \Delta \phi \Delta E\Delta t
    \nonumber \\ && + \sum_{l'=1}^L \sum _{m'=1}^M \sum _{g'=1}^G v_{g'} \sigma'
f(x_{i},y_j,z_k,\mu_{l'},\phi_{m'},E_{g'},\mu_l,\phi_m,E_g)
   n' \Delta \mu \Delta \phi  \Delta E  \Delta t \nonumber \\
&& + \sqrt{Q(x_i,y_j,z_k,E_g,t) \Delta x \Delta y \Delta z \Delta
\mu \Delta \phi \Delta E \Delta t} \; \eta^{(Q)}_{i,j,k,l,m,g}
\nonumber
\\ && - \sqrt{ v_g \sigma_c
n(x_{i},y_j,z_k,\mu_l,\phi_m,E_g,t) \Delta t} \;
\eta_{i,j,k,l,m,g}^{(c)} \nonumber \\ && - \sum_{l'=1}^L \sum
_{m'=1}^M \sum _{g'=1}^G \dfrac{1}{\hat c} \sqrt{v_g \hat \sigma
\hat f n \Delta \mu \Delta \phi \Delta E \Delta t } \;
\eta^{(tr)}_{i,j,k,l,m,g,l'm',g'} \nonumber \\&& + \sum_{l'=1}^L
\sum _{m'=1}^M \sum _{g'=1}^G \sqrt{v_{g'} \sigma' f' n' \Delta
\mu \Delta \phi
   \Delta E \Delta t} \; \eta^{(tr)}_{i,j,k,l',m',g',l,m,g}
   \nonumber
\end{eqnarray}
for $i=1,2, \dots, I$, $j=1,2, \dots, J$, $k=1,2, \dots, K$,
 $l=1,2, \dots, L$, $m=1,2, \dots, M$, and $g=1,2, \dots, G$, where
 $\eta^{(Q)}_{i,j,k,l,m,g}$, $\eta^{(c)}_{i,j,k,l,m,g}$, and
$\eta^{(tr)}_{i,j,k,l,m,g,l',m',g'}$ are independent normally
distributed numbers with zero mean and unit variance processes for
each value of $i,j,k,l,m,g,l',m',g'$. In Equation (\ref{eq2.3a}),
$\hat f= \hat
f(x_{i},y_j,z_k,\mu_l,\phi_m,E_g,\mu_{l'},\phi_{m'},E_{g'})$ and
 $f'=
f(x_{i},y_j,z_k,\mu_{l'},\phi_{m'},E_{g'},\mu_l,\phi_m,E_g).$
 For small $\Delta t$, the stochastic difference system (\ref{eq2.3a}) has
the same mean and mean square changes as the
discrete stochastic model defined by Table 1. 

Stochastic difference system (\ref{eq2.3a}) is an Euler-Maruyama
approximation to a certain It\^{o} SDE system (e.g.,
\cite{eallenbook,allenallenetal,allen03}) which has the form:
\begin{eqnarray}
 && \dfrac{d n(x_i,y_j,z_k,\mu_l,\phi_m,E_g,t)}{dt} = \label{eq2.3}
 \\ && + \mu_x v_g (n(x_{i-1},y_j,z_k,\mu_l,\phi_m,E_g,t)
  - n(x_{i},y_j,z_k,\mu_l,\phi_m,E_g,t))/  \Delta x \nonumber \\
&& + \mu_y v_g (n(x_i,y_{j-1},z_k,\mu_l,\phi_m,E_g,t)
 - n(x_{i},y_j,z_k,\mu_l,\phi_m,E_g,t))
/\Delta y \nonumber
\\ && + \mu_z v_g (n(x_i,y_j,z_{k-1},\mu_l,\phi_m,E_g,t)
 -  n(x_{i},y_j,z_k,\mu_l,\phi_m,E_g,t))
/\Delta z \nonumber \\ && + Q(x_i,y_j,z_k,E_g,t) \Delta x \Delta y
\Delta z \Delta \mu \Delta \phi \Delta E  - v_g \sigma_c
n(x_{i},y_j,z_k,\mu_l,\phi_m,E_g,t)  \nonumber \\
 && - \sum_{l'=1}^L
\sum _{m'=1}^M \sum _{g'=1}^G \dfrac{1}{\hat c} v_g \hat \sigma \hat
f(x_{i},y_j,z_k,\mu_l,\phi_m,E_g,\mu_{l'},\phi_{m'},E_{g'})
   n \Delta \mu \Delta \phi \Delta E
    \nonumber \\ && + \sum_{l'=1}^L \sum _{m'=1}^M \sum _{g'=1}^G v_{g'} \sigma'
f(x_{i},y_j,z_k,\mu_{l'},\phi_{m'},E_{g'},\mu_l,\phi_m,E_g)
   n' \Delta \mu \Delta \phi  \Delta E   \nonumber \\
&& + \sqrt{Q(x_i,y_j,z_k,E_g,t) \Delta x \Delta y \Delta z \Delta
\mu \Delta \phi \Delta E} \; \dfrac{dW^{(Q)}_{i,j,k,l,m,g}(t)}{dt}
\nonumber
\\ && - \sqrt{ v_g \sigma_c
n(x_{i},y_j,z_k,\mu_l,\phi_m,E_g,t)} \;
\dfrac{dW^{(c)}_{i,j,k,l,m,g}(t)}{dt} \nonumber \\ && -
\sum_{l'=1}^L \sum _{m'=1}^M \sum _{g'=1}^G \dfrac{1}{\hat c}
\sqrt{v_g \hat \sigma \hat f n \Delta \mu \Delta \phi \Delta E }
 \dfrac{dW^{(tr)}_{i,j,k,l,m,g,l'm',g'}(t)}{dt}  \nonumber \\&&
+ \sum_{l'=1}^L \sum _{m'=1}^M \sum _{g'=1}^G \sqrt{v_{g'} \sigma'
f' n' \Delta \mu \Delta \phi
   \Delta E} \; \dfrac{dW^{(tr)}_{i,j,k,l',m',g',l,m,g}(t)}{dt}
   \nonumber
\end{eqnarray}
for $i=1,2, \dots, I$, $j=1,2, \dots, J$, $k=1,2, \dots, K$,
 $l=1,2, \dots, L$, $m=1,2, \dots, M$, and $g=1,2, \dots, G$, where
 $W^{(Q)}_{i,j,k,l,m,g}(t)$, $W^{(c)}_{i,j,k,l,m,g}(t)$, and
$W^{(tr)}_{i,j,k,l,m,g,l',m',g'}(t)$ are independent Wiener
processes for each value of $i,j,k,l,m,g,l',m',g'$. In Equation
(\ref{eq2.3}), $\hat f= \hat
f(x_{i},y_j,z_k,\mu_l,\phi_m,E_g,\mu_{l'},\phi_{m'},E_{g'})$ and
 $f'=
f(x_{i},y_j,z_k,\mu_{l'},\phi_{m'},E_{g'},\mu_l,\phi_m,E_g).$
 For small $\Delta t$, the stochastic system (\ref{eq2.3}) has
approximately the same mean and mean square changes as the
discrete stochastic model defined by Table 1. 

Before the intervals in space, energy, and direction can be
allowed to go to zero so that the SDE system will approach an
SPDE, the Wiener processes need to be replaced with appropriate
Brownian sheets. Introduced now are multidimensional Brownian
sheets $W^{(c)}(x,y,z,\mu,\phi,E,t)$,
$W^{(Q)}(x,y,z,\mu,\phi,E,t)$, and
$W^{(tr)}(x,y,z,\mu,\phi,E,\mu',\phi',E',t)$.  For example,
$W^{(c)}(x,y,z,\mu,\phi,E,t)$ is  an independent seven-dimensional
Brownian sheet in variables $x,y,z,\mu,\phi,E,t$. The Wiener
processes in (\ref{eq2.3}) are now replaced by equivalent forms
involving Brownian sheets after which the spatial, angular, and
energy intervals will be allowed to approach zero. Specifically,
in (\ref{eq2.3}), let
\begin{eqnarray*}
&& \dfrac{dW^{(Q)}_{i,j,k,l,m,g}(t)}{dt} = \dfrac{1}{\sqrt{\Delta
x \Delta y \Delta z \Delta \mu \Delta \phi \Delta E}}
\int_{x_i}^{x_{i+1}} \int_{y_j}^{y_{j+1}} \int_{z_k}^{z_{k+1}}
\int_{\mu_l}^{\mu_{l+1}} \int_{\phi_m}^{\phi_{m+1}}  \\ &&
\int_{E_g}^{E_{g+1}} \dfrac{\partial^{7}
W^{(Q)}(x,y,z,\mu,\phi,E,t)}{\partial x
\partial y \partial z \partial \mu \partial \phi
\partial E \partial t} \; dE \,d \phi \, d  \mu \, dz \, dy \, dx
\end{eqnarray*}
\begin{eqnarray*}
&& \dfrac{dW^{(c)}_{i,j,k,l,m,g}(t)}{dt} =  \dfrac{1}{\sqrt{\Delta
x \Delta y \Delta z \Delta \mu \Delta \phi \Delta E}}
\int_{x_i}^{x_{i+1}}  \int_{y_j}^{y_{j+1}} \int_{z_k}^{z_{k+1}}
\int_{\mu_l}^{\mu_{l+1}} \int_{\phi_m}^{\phi_{m+1}} \\ &&
\int_{E_g}^{E_{g+1}} \dfrac{\partial^{7} W^{(c)}(x,y,z,\mu, \phi,
E,t)}{\partial x
\partial y
\partial z \partial \mu \partial \phi \partial E \partial t} \; dE
\, d\phi \,  d\mu \, dz \, dy \, dx
\end{eqnarray*}
\begin{eqnarray*}
&& \dfrac{dW^{(tr)}_{i,j,k,l'm',g',l,m,g}(t)}{dt} =
\dfrac{1}{\sqrt{\Delta S}} \int_{x_i}^{x_{i+1}}
\int_{y_j}^{y_{j+1}} \int_{z_k}^{z_{k+1}}
\int_{\mu_{l'}}^{\mu_{l'+1}} \int_{\phi_{m'}}^{\phi_{m'+1}}
\int_{E_{g'}}^{E_{g'+1}} \\ && \int_{\mu_l}^{\mu_{l+1}}
\int_{\phi_m}^{\phi_{m+1}} \int_{E_g}^{E_{g+1}}
\dfrac{\partial^{10} W^{(tr)}(x,y,z,\mu', \phi', E',\mu, \phi,
E,t)}{\partial x
\partial y \partial z  \partial \mu' \partial \phi' \partial E' \partial \mu
\partial \phi \partial E \partial t} \; dE \cdots dx
\end{eqnarray*}
where $\Delta S = \Delta x \Delta y  \Delta z \Delta \mu' \Delta
\phi' \Delta E' \Delta \mu \Delta \phi \Delta E$.  These
equivalent expressions are now substituted into (\ref{eq2.3}) and
$ n$ is replaced with $ N \Delta x \Delta y \Delta z \Delta \mu
\Delta \phi \Delta E$. Next, $\Delta x, \Delta y, \Delta z, \Delta
\mu, \Delta \phi,$ and $\Delta E$ are allowed to approach zero.
The result is a stochastic partial differential equation for
stochastic neutron transport:
\begin{eqnarray}
&& \dfrac{\partial N(x,y,z,\mu,\phi,E,t)}{\partial t} = -v \mu_x
\dfrac{\partial  N}{\partial x} -v \mu_y \dfrac{\partial
N}{\partial y} - v \mu_z \dfrac{\partial  N}{\partial z} \label{eq2.4} \\
&& + \,
Q(x,y,z,E,t) -v \sigma(x,y,z,E)  N(x,y,z,\mu,\phi,E) \nonumber \\
&& + \int_0^{E_{max}} \int_{-1}^1 \int_0^{2 \pi} v \sigma'
f(x,y,z,\mu',\phi',E', \mu, \phi,E) N' \; d\phi' d\mu' dE'
\nonumber \\ && + \sqrt{Q(x,y,z,E,t)} \dfrac{\partial^{7}
W^{(Q)}(x,y,z,\mu,\phi,E,t)}{\partial x
\partial y \partial z \partial \mu \partial \phi
\partial E \partial t}
 \nonumber \\ && - \sqrt{ v \sigma_c (x,y,z,E)
N(x,y,z,\mu,\phi, E,t)} \dfrac{\partial^{7} W^{(c)}(x,y,z,\mu,
\phi, E,t)}{\partial x \partial y \partial z \partial \mu \partial
\phi \partial E \partial t} \nonumber \\ && - \int_{-1}^1
\int_0^{2\pi} \int_0^{E_{max}} \dfrac{1}{\hat c} \sqrt{v \sigma f
N} \;
 \dfrac{\partial^{10} W^{(tr)}(x,y,z,\mu, \phi, E, \mu', \phi', E',t)}{\partial x
\partial y \partial z  \partial \mu \partial \phi \partial E \partial
\mu' \partial \phi' \partial E' \partial t} \; d E' d \Omega'
\nonumber \\ && + \int_{-1}^1 \int_0^{2\pi} \int_0^{E_{max}}
\sqrt{v' \sigma' f' N'} \;
 \dfrac{\partial^{10} W^{(tr)}(x,y,z,\mu', \phi', E', \mu, \phi, E,t)}{\partial x
\partial y \partial z  \partial \mu' \partial \phi' \partial E' \partial \mu
\partial \phi \partial E \partial t} \; d E' d \Omega'
 \nonumber
\end{eqnarray}
where, for convenience, $N = N(x,y,z,\mu,\phi,E,t)$,
$N'=N(x,y,z,\mu',\phi',E',t)$, $\sigma=\sigma(x,y,z,E)$,
$\sigma'=\sigma(x,y,z,E')$, $\hat \sigma(x,y,z,E) =
\sigma(x,y,z,E) - \sigma_c(x,y,z,E)$, $\hat c = \hat c(x,y,z,E)$,
$\mu_x = \mu$, $\mu_y = \sin (\cos^{-1} \mu) \cos \phi$, $\mu_z =
\sin (\cos^{-1} \mu)$ $ \sin \phi$, $d \Omega' =  d\phi' d\mu'$,
and the equation is valid for positive or negative values of the
direction cosines $\mu_x,\mu_y, \mu_z$. Notice that $N$ in
(\ref{eq2.4}) is stochastic, i.e., $N$ is not the expected neutron
angular density. That is, each solution of (\ref{eq2.4}) is one
possible realization or sample path of the neutron angular
density. There are an infinite number of these sample path
solutions, as occurs in nature, and the average of these random
solutions is equal to the expected neutron angular density $\bar
N$. Finally, notice that (\ref{eq2.4}) generalizes (\ref{eq2.1}).
If the stochastic terms are set equal to zero, then (\ref{eq2.4})
is identical to (\ref{eq2.1}). Of course, for the deterministic or
the stochastic version of the neutron transport equation, the
initial condition and boundary conditions must be specified.

As illustrated in Section \ref{calc} for two special cases of
 (\ref{eq2.4}), the stochastic transport equation (\ref{eq2.4}) can
be solved computationally by discretizing position, direction, and
energy and then approximating the resulting system of It\^{o}
stochastic differential equations in time. In effect, SDE system
(\ref{eq2.3}) is computationally solved using an appropriate
stochastic difference system such as (\ref{eq2.3a}).

\section{Special Cases of the Stochastic Neutron Transport Equation}
\label{cases}

   In this section, for illustrative purposes, two special cases
   of (\ref{eq2.4}) are considered. First, a
   one-dimensional stochastic neutron transport equation with
   isotropic scattering is described. Second, a stochastic
   partial differential equation is given for neutrons interacting
   in a homogeneous medium.

      Consider the parallelepiped  region $0 \le x \le x_{max}$,
      $0 \le y \le y_{max}$,
 and $0 \le z \le z_{max}$ where the medium is uniform with respect
 to the spatial variables $y$ and $z$. Assume that, except for the
 left face and right face defined
 by $x=0$ and  $x=x_{max}$, respectively, the neutrons are reflected back at
 the other faces,
 i.e., there are reflecting boundary conditions at the four faces
 except for the left
 and right faces. Equation (\ref{eq2.4}) is integrated over
 $0 \le y \le y_{max}$, $0 \le z \le z_{max}$, and angle
 $\phi$.  Then, the equation reduces to
\begin{eqnarray}
&& \dfrac{\partial N(x,\mu,E,t)}{\partial t} = -v \mu
\dfrac{\partial  N}{\partial x}
 + Q(x,E,t) -v \sigma(x,E)  N(x,\mu,E,t) \label{eq3.1} \\&&
   + \int_0^{E_{max}} \int_{-1}^1  v'
\sigma(x,E') f(x,\mu',E',\mu,E) N(x,\mu',E',t) \;  d\mu' dE'
\nonumber  \\ && + \sqrt{Q(x,E,t)} \dfrac{\partial^{4}
W^{(Q)}(x,\mu,E,t)}{\partial x \partial \mu
\partial E \partial t} - \sqrt{ v \sigma_c (x,E) N} \dfrac{\partial^{4}
W^{(c)}(x,\mu, E,t)}{\partial x
\partial \mu  \partial E \partial t} \nonumber \\ && - \int_{-1}^1
\int_0^{E_{max}} \dfrac{1}{\hat c(x,E)}\sqrt{v \sigma f N} \;
 \dfrac{\partial^{6} W^{(tr)}(x,\mu, E, \mu', E',t)}{\partial x
 \partial \mu  \partial E \partial \mu' \partial E' \partial t} \; d E'
  d\mu' \nonumber \\ && +
\int_{-1}^1 \int_0^{E_{max}} \sqrt{v' \sigma' f' N'} \;
\dfrac{\partial^{6} W^{(tr)}(x,\mu', E', \mu, E,t)}{\partial x
\partial \mu'  \partial E' \partial \mu
 \partial E \partial t} \; d E' d\mu'
 \nonumber
\end{eqnarray}
where $N(x,\mu,E,t)$ is the number of neutrons per unit length per
unit angle per unit energy at position $x$ with direction $\mu$ and
energy $E$ at time $t$, $\hat c(x,E) = \int_0^{E_{max}} \int_{-1}^1
 \hat f(x,\mu,E,\mu',E') d \mu' d E'$ is the mean number of neutrons
emerging per non-capture collision of neutrons of energy $E$ at
position $x$, $N' =N(x, \mu', E',t)$, and $\sigma' =
\sigma(x,E')$. In Equation (\ref{eq3.1}), $ f'=
f(x,\mu',E',\mu,E)$ and $f= f(x,\mu,E,\mu',E')$. Equation
(\ref{eq3.1}) is a stochastic neutron transport equation for
one-dimensional slab geometry.

Furthermore, assuming a single energy, isotropic scattering, and
only capture and scattering interactions, the above equation
becomes:
\begin{eqnarray}
&& \dfrac{\partial N(x,\mu,t)}{\partial t} = -v \mu
\dfrac{\partial N}{\partial x}
 + Q -v \sigma(x)  N
 +  \dfrac{1}{2} \int_{-1}^1  v
\sigma_s(x) N(x,\mu',t) \;  d \mu' \label{eq3.2}
\\ &&   + \sqrt{Q(x,t)}
\dfrac{\partial^{3} W^{(Q)}(x,\mu,t)}{\partial x \partial \mu
\partial t} + \int_{-1}^1 \sqrt{\frac{1}{2} v \sigma_s(x) N(x, \mu',t)} \;
 \dfrac{\partial^{4} W^{(tr)}(x,\mu', \mu,t)}{\partial x
 \partial \mu'  \partial \mu \partial t} \; d\mu' \nonumber \\
 && - \sqrt{ v \sigma_c (x)
N} \dfrac{\partial^{3} W^{(c)}(x,\mu,t)}{\partial x
\partial \mu  \partial t}  - \int_{-1}^1
\sqrt{\frac{1}{2} v \sigma_s(x) N} \;
 \dfrac{\partial^{4} W^{(tr)}(x,\mu, \mu', t)}{\partial x
 \partial \mu \partial \mu' \partial t} \;  d\mu'.
 \nonumber
\end{eqnarray}
where $N=N(x,\mu,t)$, $Q=Q(x,t)$, and $\sigma_s(x)$ is the
scattering cross section at position $x$. Equation (\ref{eq3.2})
is a stochastic neutron transport equation for mono-energetic
transport in one-dimensional plane geometry with isotropic
scattering.

Consider again the parallelepiped  region $0 \le x \le x_{max}$,
     $0 \le y \le y_{max}$, and $0 \le z \le z_{max}$ where the medium
     is uniform with respect
 to all the spatial variables $x$, $y$, and $z$. Assume that the neutrons are
 reflected back at all six faces,
 i.e., there are reflecting boundary conditions at all the faces.
 Equation (\ref{eq2.4}) is integrated over the
 volume $0 \le x \le x_{max}$, $0 \le y \le y_{max}$, and $0 \le z \le z_{max}$,
 and over the angles
 $\phi$ and $\mu$.  Then,
   equation (\ref{eq2.4}) reduces to
\begin{eqnarray}
&& \dfrac{\partial N(E,t)}{\partial t} = -v \sigma(E) N
 + Q(E,t) + \sqrt{Q(E,t)} \dfrac{\partial^{2}
W^{(Q)}(E,t)}{\partial E \partial t}  \label{eq3.3}
\\ && + \int_0^{E_{max}}  v'
\sigma(E') f(E',E) N(E',t) \;  d E' - \sqrt{ v \sigma_c (E) N}
 \dfrac{\partial^{2} W^{(c)}(E,t)}{\partial E \partial t} \nonumber \\ &&
 -
\int_0^{E_{max}}\dfrac{1}{\hat c(E)} \sqrt{v \sigma(E) f(E,E')
N(E,t)} \;
 \dfrac{\partial^{3} W^{(tr)}(E, E',t)}{\partial E  \partial E'
 \partial t} \; d E'  \nonumber \\ &&  +
 \int_0^{E_{max}} \sqrt{v' \sigma(E') f(E',E) N(E',t)} \;
 \dfrac{\partial^{3} W^{(tr)}(E', E,t)}{\partial E' \partial E \partial t}
  \; d E'
 \nonumber
\end{eqnarray}
where $N=N(E,t)$ is the number of neutrons per unit energy and
$Q(E,t)$ is equal to the number of source neutrons per unit energy
per unit time.

In a homogeneous medium, with only capture and scattering
interactions, the above equation becomes:
\begin{eqnarray}
&& \dfrac{\partial N(E,t)}{\partial t} = -v \sigma(E) N
 + Q(E,t) + \sqrt{Q(E,t)} \; \dfrac{\partial^{2}
W^{(Q)}(E,t)}{\partial E \partial t}  \label{eq3.4}
\\ && +\int_0^{E_{max}} v' \sigma(E') f(E',E) N(E',t) \;  dE'
-\sqrt{ v \sigma_c (E) N} \; \dfrac{\partial^{2} W^{(c)}(
E,t)}{\partial E \partial t} \nonumber \\ && - \int_0^{E_{max}}
\sqrt{v \sigma(E) f(E,E') N(E,t)} \;
 \dfrac{\partial^{3} W^{(tr)}(E, E',t)}{\partial E  \partial E'
 \partial t} \; d E' \nonumber \\ &&  +
\int_0^{E_{max}} \sqrt{v' \sigma(E') f(E',E) N(E',t)} \;
 \dfrac{\partial^{3} W^{(tr)}(E', E, t)}{\partial E' \partial E \partial t}
  \; d E'
 \nonumber
\end{eqnarray}
where $ v' \sigma(E') f(E',E)$ is the probability of a neutron
scattering from energy $E'$ to $E$ per unit energy per unit time.
Equation (\ref{eq3.4}) is a stochastic neutron transport equation
for a homogeneous medium with only capture and scattering
collisions.

\section{Comparison With Monte Carlo Calculations}
\label{calc}

   In this section, the stochastic difference equations
derived in the previous sections for neutron transport are
numerically solved and compared with independent Monte Carlo
computations. Two cases are considered.
 First, a problem is studied involving mono-energetic
neutron transport in a slab with captures and isotropic scatters.
The stochastic transport equation in this case is given by
equation (\ref{eq3.2}). Second, a  homogeneous medium is studied
where the neutrons experience captures or scatters with energy
changes. The stochastic transport equation for this problem is
(\ref{eq3.4}).

    In the first problem, it is assumed that
1000 neutrons per second begin entering the left side of a slab of
width one unit at time $t=0.0$.  The velocity of the neutrons is
$v=0.1$. The slab is homogeneous and the scattering and capture
cross sections in the slab are assumed to be $\sigma_s=5.0$ and
$\sigma_c=0.10$. Therefore, for this problem, the slab has width
$1/\sigma_c=10$ absorption mean free paths. (As slab width
increases, the leakage decreases but it is likely that the
coefficient of variation in the leakage increases.) The neutrons
isotropically enter the slab on the left side, $x=0$, from time
$t=0.0$ until $t=50.0$. After time $t=50.0$, the neutrons no
longer enter the slab and the neutrons in the slab eventually are
absorbed or escape. There is no neutron source for this problem.
The time dependence of the exiting fluxes on the left and right
sides of the slab are of interest in this problem. To study this
problem computationally, equation (\ref{eq3.2}) needs to be solved
numerically.

To define a numerical method for equation (\ref{eq3.2}), the
interval $[0,1]$ in $x$ is divided into $I$ intervals
 $[x_{i},x_{i+1}]$, for $i=1,2, \dots, I$ where $x_i=(i-1) \Delta x$,
 and $\Delta x =1/I$. In addition,  the interval $[-1,1]$ in direction $\mu$ is
 divided into $J$ equal intervals of width $\Delta \mu = 2/J$
  and time is discretized where
 $t_k = k \Delta t, k=0, 1, 2, \dots$. Considering (\ref{eq3.2})
 at position $x_i$ and direction $\mu_j$
 and using an explicit approximation in time $t$ along with an
 upwind differencing approach suggests the numerical procedure:
     \begin{eqnarray}
&&  n_{i,j,k+1} = n_{i,j,k} + \left\{\begin{array}{ll}  \mu_j v
n_{i-1,j,k}\Delta t/\Delta x -
\mu _j v n_{i,j,k}\Delta t/\Delta x, \;\; \text{for} \;\; \mu_j >0 \\
 \mu_j v n_{i,j,k}\Delta t/\Delta x -
\mu _j v n_{i+1,j,k}\Delta t/\Delta x, \;\; \text{for} \;\; \mu_j
<0
\end{array}\right. \label{eq4.1}
 \\ &&  - v \sigma_{i}
n_{i,j,k} \Delta t + \sum_{m=1}^J \dfrac{1}{2} \sigma_{s,i} v
n_{i,m,k} \Delta \mu \Delta t - \sqrt{v \sigma_{c,i} n_{i,j,k}
\Delta t}\; \eta_{i,j,k}^{(c)} \nonumber
\\ && - \sum_{m=1}^J \sqrt{ \dfrac{1}{2} \sigma_{s,i} v n_{i,j,k}
\Delta \mu \Delta t} \; \eta_{i,j,m,k}^{(tr)} + \sum_{m=1}^J
\sqrt{\dfrac{1}{2} \sigma_{s,i} v n_{i,m,k} \Delta \mu \Delta t}
\; \eta_{i,m,j,k}^{(tr)} \nonumber
\end{eqnarray}
where   $n_{i,j,k} \approx n(x_i,\mu_j, t_k)$ is the number of
neutrons at position $x_i$ in direction $\mu_j$ at time $t_k$ and
$\sigma_{s,i} = \sigma_s(x_i)$. Also, $\eta_{i,j,m,k}^{(tr)}$ and
$\eta_{i,j,k}^{(c)}$ are independent Gaussian $\mathcal
N(0,1)$-distributed numbers for each $i,j,k,m$. Notice that
(\ref{eq4.1}) is an Euler-Maruyama approximation
\cite{gard,kloedenplaten,kps} to the system of It\^{o}
differential equations (\ref{eq2.3}) and is a special case of the
stochastic difference system (\ref{eq2.3a}).

\begin{figure}[ht]
$$\includegraphics[height=2.2in,width=4in]{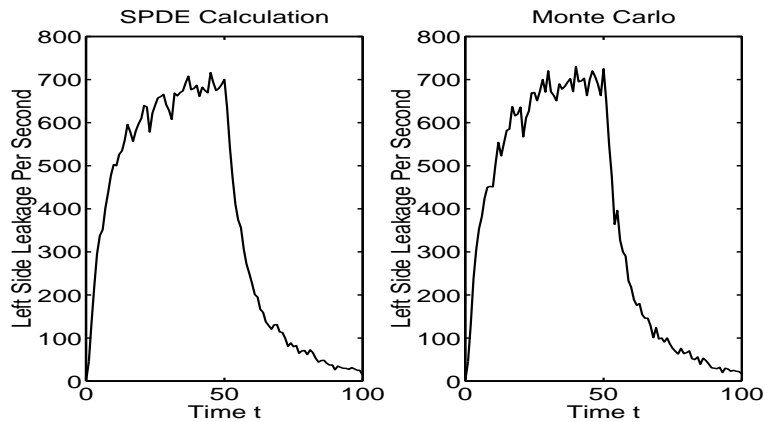}$$
\caption{Calculated left leakages per second from time $t=0$ to
$t=100$ for one sample path using Monte Carlo and SPDE
(\ref{eq3.2}).}\label{figslab1}
\end{figure}

\begin{figure}[ht]
$$\includegraphics[height=2.2in,width=4in]{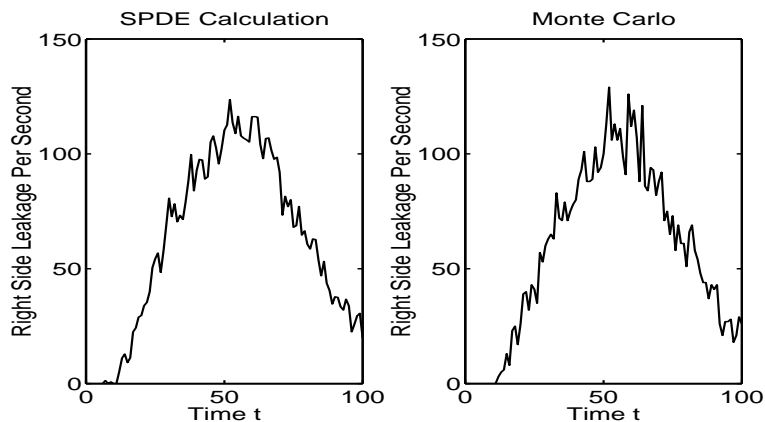}$$
\caption{Calculated right leakages per second from time $t=0$ to
$t=100$ for one sample path using Monte Carlo and SPDE
(\ref{eq3.2}).}\label{figslab2}
\end{figure}

     The problem is solved numerically using two independent computational
procedures, i.e, numerical solution of the SPDE is compared with
Monte Carlo calculations. Equation (\ref{eq4.1}) is solved
computationally with $I=80$ equal intervals in position $x$ and
$J=40$ equal angular intervals. The value chosen for the time
interval is $\Delta t=0.125$. In the Monte Carlo calculations,
1000 neutrons per second enter isotropically on the left side.
Each neutron is followed individually in the Monte Carlo procedure
with each neutron checked for a scatter, a leakage, or an
absorption at each time step of 0.1 seconds. Calculational results
for 100 sample paths using the two independent computational
approaches are given in Table 2. 
The means and standard deviations of the calculated number of
neutrons escaping from the left side and from the right side are
given for the time interval $t=49$ to $t=50$. The two approaches
agree well. In Figures \ref{figslab1} and \ref{figslab2}, the
calculated leakages for one sample path are compared for the two
approaches.

     The second problem involves neutrons slowing down in a
     homogeneous medium. It is assumed that at time $t=0.0$,
there are 0 neutrons with energy between 0 eV  and 10 eV and 400
neutrons with energy between 10 eV and 20 eV. Also, there is
present a constant source, $Q(E)$,  of neutrons where $Q(E) =
\left\{\begin{array}{ll} 22 \;\;
 \text{for} \;\; 10 \le E \le 20 \\ 0
\;\; \text{for} \;\; 0 \le E \le 10.\nonumber
\end{array}\right.$ Thus, $\int_0^{20} Q(E) \, d E = 220$ source
neutrons are produced per second. For this problem, the source,
$Q(E)$, is assumed to be constant and does not vary randomly. It
is furthermore assumed that the total cross section satisfies
$v(E) \sigma(E) = 1.0 \;\; \text{for} \;\; 0 \le E \le 20$, the
capture cross section has the form $ v(E) \sigma_c(E) =
\left\{\begin{array}{ll} 0.1 \;\;
 \text{for} \;\; 10 \le E \le 20 \\ 1.0
\;\; \text{for} \;\; 0 \le E \le 10
\end{array}\right.$ and the
scattering cross section has the form $ v(E') \sigma(E') f(E',E) =
\left\{\begin{array}{ll} 0.045 \;\;
 \text{for} \;\;  10 \le E' \le 20 \\ 0.0
\;\; \text{for} \;\; 0 \le E' \le 10.\nonumber
\end{array}\right.$  For convenience, the product of the speed, $v(E')$,
with the cross sections are given. The scattering kernel,
$\sigma(E') f(E',E)$, is thus assumed to be piecewise continuous
and proportional to $1/\sqrt{E'}$ rather than, for example, to
$1/E'$ such as for a hydrogen-moderated system. Also, notice that
the cross sections are consistent in the sense that
$$ v(E') \sigma(E') = v(E') \sigma_c(E') + \int_0^{20} v(E') \sigma(E',E)
f(E',E)\, dE = 1 \;\;\; \text{for all} \;\;\; 0 \le E' \le 20.$$
Furthermore, for this problem, it can be shown that the mean
number of neutrons with energies between 10 eV and 20 eV is equal
to 400 for time $t \ge 0$ and the mean number of low-energy
neutrons with energies between 0 eV and 10 eV
 approaches 180 as time $t$ increases. In this problem, the stochastic behavior
 of the
 number of neutrons with energies  between 0 eV and 10 eV is of
 interest.

     This problem is solved numerically using two independent computational
procedures. Specifically, numerical solution of  SPDE
(\ref{eq3.4}) is compared with Monte Carlo calculations. In the
numerical solution of Equation (\ref{eq3.4}), $G=20$ energy groups
of equal width $\Delta E = 20/G$ are used. Equation (\ref{eq3.4})
is solved using the following stochastic difference equations at
discrete times $t_k = k \Delta t$ where $\Delta t = 0.02$:
    \begin{eqnarray}
&&  n_{g,k+1} = n_{g,k}
 +  q_{g,k} \Delta E \Delta t - v_g \sigma_{g}
n_{g,k} \Delta t + \sum_{g'=1}^G v_{g'} \sigma_{g'} f_{g',g}
n_{g',k} \Delta E \Delta t  \label{eq4.2}
\\ && - \sqrt{v_g \sigma_{c,g} n_{g,k} \Delta
t}\; \eta_{g,k}^{(c)}  - \sum_{g'=1}^G \sqrt{ v_g  \sigma_{g}
f_{g,g'} n_{g,k} \Delta E \Delta t} \; \eta_{g,g',k}^{(tr)}
\nonumber \\ && + \sum_{g'=1}^G \sqrt{ v_{g'} \sigma_{g'} f_{g',g}
n_{g',k} \Delta E \Delta t} \; \eta_{g',g,k}^{(tr)} \nonumber
\end{eqnarray}
where $n_{g,k} \approx n(E_g,t_k)$ is the number of neutrons in
the $g$th energy group at time $t_k = k \Delta t$  and $\sigma_g
f_{g,g'} = \sigma(E_g) f(E_g, E_{g'})$. Also, $\eta_{g,k}^{(c)},
\eta_{g,g',k}^{(tr)} \sim \mathcal N(0,1)$ are independent
normally distributed numbers with mean 0 and variance 1 for each
$g,g',k$. Difference system (\ref{eq4.2}) is a special case of the
stochastic difference system (\ref{eq2.3a}). In the Monte Carlo
procedure, the neutron population in each energy group is checked
at each time step for an absorption, an energy group change, or
for an addition from the neutron source. This procedure is
continued for each time step until the final time $t=2$.

\begin{figure}[ht]
$$\includegraphics[height=2.2in,width=4in]{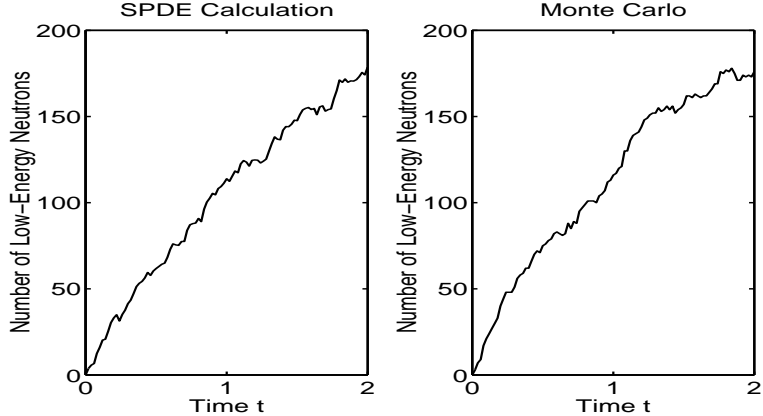}$$
\caption{Calculated number of neutrons with energy less than 10.0
from time $t=0$ to $t=2$ for one sample path using Monte Carlo and
SPDE (\ref{eq3.4}).} \label{figenergy1}
\end{figure}

 Calculational results
for 100 sample paths using the two independent computational
approaches are given in Table 3.   
The means and standard deviations of the calculated number of
neutrons with energies between 0 eV and 10 eV and with energies
between 10 eV and 20 eV are given for time $t=2$.  The two
different computational approaches agree well. In Fig.
\ref{figenergy1}, the calculated number of neutrons with energies
between 0 eV and 10 eV are given from time $t=0$ to $t=2$ for one
sample path for each calculational method. Again, the results are
very similar for the two different calculational procedures.

\section{Applications}
\label{app}

   Three  possible applications of the stochastic
 neutron transport equations (\ref{eq2.4}) are
 in development of new computational methods for solving
 neutron transport problems, in sensitivity or
 perturbation analysis, and in testing of computational and analytical methods.
 Numerical methodology for stochastic partial differential equations is advancing rapidly. It is
 probable that numerical techniques for solving stochastic
 neutron transport equations may eventually out-compete Monte Carlo
 techniques. Stochastic neutron transport equations
 can also be useful, for example, in perturbation studies to
 evaluate the effects of small changes in reactor parameters. In addition,
 solutions of stochastic neutron transport equations can provide independent
 checks on numerical or analytical approaches such as Monte Carlo
 methods.

      Consider a very simple example of applying
stochastic neutron transport equations in a perturbation study.
Consider energy-dependent transport in a homogeneous medium with
the neutron density given by (\ref{eq3.4}). Assume, for this
example, that there is no neutron source and there are no
collisions other than capture collisions where the capture cross
section varies with time. The stochastic transport equation, for
this problem, reduces to:
\begin{eqnarray}
&& \dfrac{\partial N(E,t)}{\partial t} = -v \sigma_c(E,t) N(E,t)
-\sqrt{ v \sigma_c (E,t) N(E,t)} \; \dfrac{\partial^{2} W(
E,t)}{\partial E
\partial t}.  \label{eq3.5}
\end{eqnarray}
Suppose that $\sigma_c(E,t)$ is perturbed to  $\tilde
\sigma_c(E,t) = \sigma_c(E,t) + \Delta \sigma_c(E,t)$ for $t \ge
0$. We wish to estimate, for the perturbation, the total number of
neutrons as well as the change in the variability in this number
with time $t$.
 In particular, if the perturbed neutron density is $ N(E,t)
 + \Delta N(E,t)$, we wish to estimate the change in the total number of neutrons as a function
 of time, i.e.,  $\Delta n(t) = \int_0^{\infty} \Delta N(E,t) \, dE$.
 To derive equations for this quantity, let $E_k =k \Delta E$ where $\Delta E$
 is a small energy interval. Integrating (\ref{eq3.5}) over
 the interval $[E_k, E_{k+1}]$, then
 \begin{equation}
  \dfrac{d n_k (t)}{d t} = -v_k \sigma_c(E_k,t) n_k(t) -
 \sqrt{v_k \sigma_c(E_k,t) n_k(t)} \, \dfrac {d W_k(t)}{d t} \label{eq5.1}
 \end{equation}
 for $k=0, 1, \dots$ where $n_k(t) = \int_{E_k}^{E_{k+1}} N(E,t) \, dE$.
 In addition, using It\^{o}'s formula (e.g., \cite{gard}),
\begin{equation}
  \dfrac{d n_k^2 (t)}{d t} = -2v_k \sigma_c(E_k,t) n_k^2(t) + v_k \sigma_c(E_k,t)
  n_k(t) - 2 n_k(t) \sqrt{v_k \sigma_c(E_k,t) n_k(t)} \, \dfrac {d W_k(t)}{d t}. \label{eq5.2}
 \end{equation}
 From (\ref{eq5.1}) and (\ref{eq5.2}), expressions for $E(n_k(t))$
 and $E(n_k^2(t))$ are readily obtained and then, $E(n(t)) =  E(\sum_k n_k(t))$ and
 $E(n^2(t)) = E((\sum_k n_k(t))^2)$.  Finally, expressions for
 $ E(n(t) + \Delta n(t)) - E(n(t))$ and
 $ \text{Var}(n(t) + \Delta n(t)) - \text{Var}(n(t)) $ are derived as:
 \begin{equation}
  E(n(t) + \Delta n(t)) - E(n(t)) = \int_0^{\infty} N(E,0)\left(e^{-\int_0^t v \tilde \sigma_c(E,s) ds} -
    e^{-\int_0^t v \sigma_c(E,s) ds}\right) \, dE \label{eq5.3}
 \end{equation}
 and
\begin{eqnarray}
 && \text{Var}(n(t) + \Delta n(t)) - \text{Var}(n(t))  = \nonumber
 \\ &&
 \int_0^{\infty} N(E,0) e^{-2\int_0^t v \tilde \sigma_c(E,s) ds}
 \left(\int_0^t v \tilde \sigma_c(E,s)
 e^{\int_0^s v \tilde \sigma_c(E,r) dr} \, ds \right)\, dE  \nonumber \\
 && -
 \int_0^{\infty} N(E,0) e^{-2\int_0^t v \sigma_c(E,s) ds}
 \left(\int_0^t v \sigma_c(E,s)
 e^{\int_0^s v \sigma_c(E,r) dr} \, ds \,\right) dE.
 \label{eq5.4}
 \end{eqnarray}
 Therefore, using the stochastic neutron transport equation to analyze the effect of
 the perturbation for this transport process, not only are equations derived for estimating
 the average effect of the perturbation but equations are also
 obtained for estimating the change in the variability for the
 perturbation. Indeed, for this perturbation problem, equations (\ref{eq5.3})
 and (\ref{eq5.4}) clearly indicate that the mean change and the
 change in the variability are both proportional to the initial
 number of neutrons.

\section{Conclusions and Future Directions}
\label{sum}

Stochastic difference and partial differential equations (SPDEs)
are becoming increasingly important in applied mathematics
\cite{daprato,gunzburger,holdenetal,schurz}. In the present
investigation, stochastic difference equations and an SPDE are
derived for neutron transport in a general three-dimensional
medium. In the derivation procedure, the deterministic and
stochastic terms in the differential equation system are
simultaneously derived. First, a stochastic difference system is
constructed. Next, an SDE system is derived. Finally, a particular
SPDE follows from the SDE system. The stochastic difference
equations and the SPDE for the neutron angular densities are given
by (\ref{eq2.3a}) and (\ref{eq2.4}), respectively. SPDEs for
special cases of this equation are given by (\ref{eq3.2}) and
(\ref{eq3.4}) for one-dimensional plane geometry and for a
homogeneous medium, respectively.  The stochastic equations
generalize the deterministic neutron transport equations and
include random influences due to interactions and sources. Hence,
certain random phenomena, such as fluctuations during reactor
startup, can be studied using these SPDEs. The stochastic
difference equations for neutron transport are solved numerically
and compared with independently formulated Monte Carlo methods.
The computational results between the two different numerical
methods are in good agreement supporting the accuracy of the
stochastic neutron transport derivation procedure.

 Future work may include
appropriately extending the derivations of the present
investigation to include the random influence of prompt and
delayed neutrons \cite{bellglasstone,hayesallen,hetrick}. In
addition, stochastic difference equations and an SPDE may be
developed to model the random behavior of neutron transport in
spherical and cylindrical geometries.

\section*{Acknowledgement}
This work was partially supported by NSF grant DMS-0718302.

\newpage

\begin{center}
TABLE 1 \\
Possible Changes in the Number $n =
n(x_i,y_j,z_k,\mu_l,\phi_m,E_g,t_p)$ for Time $\Delta t$
\end{center}
\begin{center}
 \begin{table}[ht]
 \begin{center}
\begin{tabular}{|c|c|c|}
 \hline
 Change $(\Delta n)$& Description & Probability in Time $\Delta t$ \\\hline
 $\mu_x v_g  n(x_{i-1},y_j,z_k,\mu_l,\phi_m,E_g,t_p)\Delta t / \Delta x$
   & In a yz-face $(\mu_x >0)$ & $1$ \\
 $-\mu_x v_g n(x_{i+1},y_j,z_k,\mu_l,\phi_m,E_g,t_p)\Delta t / \Delta x$
  & In a yz-face $(\mu_x <0)$ & $1$ \\
 $-|\mu_x| v_g n(x_{i},y_j,z_k,\mu_l,\phi_m,E_g,t_p)\Delta t / \Delta x$
 & Out a yz-face & $1$ \\
 $\mu_y v_g n(x_i,y_{j-1},z_k,\mu_l,\phi_m,E_g,t_p)\Delta t / \Delta y$
   & In a xz-face $(\mu_y >0)$ & $1$ \\
$-\mu_y v_g n(x_i,y_{j+1},z_k,\mu_l,\phi_m,E_g,t_p)\Delta t /
\Delta y$ & In a xz-face $(\mu_y <0)$ & $1$ \\
$-|\mu_y| v_g n(x_{i},y_j,z_k,\mu_l,\phi_m,E_g,t_p)\Delta t /
\Delta y$ & Out a xz-face & $1$ \\
$\mu_z v_g n(x_i,y_j,z_{k-1},\mu_l,\phi_m,E_g,t_p)\Delta t /
\Delta z$ & In a xy-face $(\mu_z >0)$ & $1$ \\
 $-\mu_z v_g n(x_i,y_j,z_{k+1},\mu_l,\phi_m,E_g,t_p)\Delta t / \Delta z$
 & In a xy-face $(\mu_z <0)$ & $1$ \\
 $-|\mu_z| v_g n(x_{i},y_j,z_k,\mu_l,\phi_m,E_g,t_p)\Delta t / \Delta z$
  & Out a xy-face & $1$ \\
 $-1 $ & Capture & $p_c$\\
 $-1/\hat c$  & Transfer out  &
 $p_{tr1}$ \\
 $1$  & Transfer in  &
 $ p_{tr2}$ \\
   $1$  & Source  &
   $p_Q $  \\
\hline
\end{tabular} \label{table1}
\end{center}
\end{table}
\end{center}

\newpage

\begin{center}
TABLE 2 \\
Monte Carlo (MC) and SPDE Calculational Results for 100 Sample \\
Paths for the Leakage for the Time Interval $t=49$ to $t=50$
\end{center}
\begin{center}
\begin{table}[ht]
\begin{center}
\begin{tabular}{|c|c|c|c|}
\hline
 Average Number & Standard & Average Number & Standard \\
Out Left Side & Deviation & Out Right Side & Deviation \\
\hline 704.93 (MC) \;\; &   23.01 (MC) \;\; &
 100.22 (MC) \;\; & 9.93 (MC) \;\; \\
 694.32 (SPDE) &  21.05 (SPDE) & 106.75 (SPDE) & 7.57 (SPDE) \\
\hline
\end{tabular}
 \label{table2}
 \end{center}
\end{table}
\end{center}

\newpage

\begin{center}
TABLE 3 \\
\noindent Monte Carlo (MC) and SPDE Results for 100 Sample Paths
for the Number of \\ Neutrons at Time $t=2.0$ With Energy Less
Than 10.0 or Between 10.0 and 20.0
\end{center}
\begin{center}
\begin{table}[ht]
\begin{center}
\begin{tabular}{|c|c|c|c|}
\hline
 Average Number & Standard & Average Number & Standard\\
With Energy & Deviation & With Energy & Deviation \\
Less Than 10.0 & & Between 10.0 and 20.0 &  \\
\hline 156.83 (MC) \;\; &   11.44 (MC) \;\; &
 399.97 (MC) \;\; & 14.39 (MC) \;\; \\
 156.97 (SPDE) &  10.28 (SPDE) & 400.92 (SPDE) & 13.42 (SPDE) \\
\hline
\end{tabular}
 \label{table3}
\end{center}
\end{table}
\end{center}

\newpage


\vspace{.5in}

\noindent {\bf Figure Captions}

\vspace{.5in}

\noindent Fig. 1. A Brownian sheet on $[0,5] \times [0,5]$.

\vspace{.3in}

\noindent Fig. 2. Calculated left leakages per second from time $t=0$ to
$t=100$ for one sample path using Monte Carlo and SPDE
(\ref{eq3.2}).

\vspace{.3in}

\noindent Fig. 3. Calculated right leakages per second from time $t=0$ to
$t=100$ for one sample path using Monte Carlo and SPDE
(\ref{eq3.2}).

\vspace{.3in}

\noindent Fig. 4. Calculated number of neutrons with energy less than 10.0
from time $t=0$ to $t=2$ for one sample path using Monte Carlo and
SPDE (\ref{eq3.4}).

\vspace{.3in}

\newpage






\begin{thebibliography}{00}


\bibitem{eallenbook} E. J. Allen,  {\it Modeling With It\^{o} Stochastic
Differential Equations}, Springer, Dordrecht, 2007.

\bibitem{eallen0708} E. J. Allen, {\it Derivation of Stochastic
Partial Differential Equations}, Stoch. Anal. Appl. 26(2008), pp.
357-378.

\bibitem{eallenjbd} E. J. Allen,  {\it Derivation of Stochastic
Partial Differential Equations for Size- and Age-Structured
Populations},  J. Bio. Dyn. 3(2009), pp. 73-86.

\bibitem{allenallenetal} E. J. Allen, L. J. S. Allen, A. Arciniega,
P. E. Greenwood, {\it Construction of equivalent stochastic
differential equation models}, Stoch. Anal. Appl. 26(2008), pp.
274-297.

\bibitem{allennovoselzhang} E. J. Allen, S. J. Novosel, Z. Zhang,
{\it Finite element and difference approximation of some linear
stochastic partial differential equations}, Stochastics and
Stochastics Reports 64(1998), pp. 117-142.

\bibitem{allen03}  L. J. S. Allen, {\it An Introduction to Stochastic Processes
with Applications to Biology}, Pearson Education Inc., Upper
Saddle River, New Jersey, 2003.


\bibitem{bella} G. I. Bell, {\it Probability distribution of neutrons
and precursors in a multiplying assembly}, Ann. Phys. 21(1963),
pp. 243-283.



\bibitem{bellb} G. I. Bell, {\it On the stochastic theory of neutron
transport}, Nucl. Sci. and Eng. 21(1965), pp. 390-401.

\bibitem{bellglasstone} G. I. Bell, S. Glasstone,  {\it Neutron Transport
Theory}, Van Nostrand Reinhold Company, New York, 1970.


\bibitem{cabana} E. M. Caba\~na, {\it The vibrating string forced by white
noise}, Z. Wahrscheinlichkeit. 15(1970), pp. 111-130.

\bibitem{daprato} G. Da Prato, L. Tubaro (Eds.), {\it Stochastic Partial
Differential Equations and Applications - VII},  CRC Press, Taylor
\& Francis Group, Boca Raton, Florida, 2006.

\bibitem{duderstadtmartin} J. Duderstadt, W. Martin, {\it Transport Theory},
John Wiley and Sons, New York, 1979.

\bibitem{gard} T. C. Gard,  {\it Introduction to Stochastic
Differential Equations}, Marcel Decker, New York, 1987.

\bibitem{gillespie} D. T. Gillespie, {\it The chemical Langevin
equation}, J. Chem. Phys. 113(2000), pp. 297-306.

\bibitem{gunzburger} M. Gunzburger, {\it Numerical methods for
stochastic PDEs}, SIAM News 40(2007), pg. 3.

\bibitem{hayesallen} J. G. Hayes and E. J. Allen, {\it Stochastic point-kinetics
equations in nuclear reactor dynamics}, Ann. Nucl. Eng. 32(2005),
pp. 572-587.

\bibitem{hetrick} D. L. Hetrick, {\it Dynamics of Nuclear Reactors},
 The University of Chicago Press, Chicago, 1971.

\bibitem{holdenetal} H. Holden, B. \O ksendal, J. Ub\o e, T. Zhang,
{\it Stochastic Partial Differential Equations: A Modeling, White
Noise Functional Approach}, Birh\"{a}user, Boston, Massachusetts,
1996.

\bibitem{palb} Y. Kitamura, L. P\'{a}l, I. P\'{a}zit, A. Yamamoto,
 Y. Yamane, {\it Some properties of zero power noise in a time-varying medium
with delayed neutrons}, Ann. Nucl. Eng. 35(2008), pp. 1621-1627.

\bibitem{kloedenplaten}  P. E. Kloeden, E. Platen, {\it Numerical
Solution of Stochastic Differential Equations}, Springer-Verlag,
New York, 1992.

\bibitem{kps} P. E. Kloeden, E. Platen, H. Schurz,
{\it Numerical Solution of SDE Through Computer Experiments},
Springer, Berlin, 1994.

\bibitem{lewismiller} E. E. Lewis, W. F. Miller, {\it Computational
Methods of Neutron Transport}, John Wiley, New York, 1984.


\bibitem{pala} L. P\'{a}l, I. P\'{a}zit, {\it Theory of neutron noise in a
temporally fluctuating multiplying medium}, Nucl. Sci. Eng.
155(2007), pp. 425-440.


\bibitem{schurz} H. Schurz, {\it Nonlinear stochastic wave equations
in $\mathbb{R}^1$ with power-law nonlinearity and additive
space-time noise}, Contemp. Math. 440(2007), pp. 223-242.

\bibitem{sharpallen} W. D. Sharp and E. J. Allen, {\it Stochastic neutron
   transport equations for rod and plane geometries},
   Ann. Nucl. Eng. 27(2000), 99-116.

\bibitem{vankampen}  N. G. van Kampen, {\it Stochastic Processes
in Physics and Chemistry}, Elsevier Science B. V., Amsterdam, The
Netherlands, 1992.

\bibitem{walsh} J. B. Walsh, {\it An Introduction to Stochastic Partial
Differential Equations}, in {\it Lecture Notes in Mathematics,
Vol. 1180}, A. Dold, B. Eckmann, eds., Springer-Verlag, Berlin,
1986, pp. 265-439.

\bibitem{williamslarsen} M. M. R. Williams, E. W. Larsen,
{\it Neutron transport in spatially random media: eigenvalue
problems}, Nucl. Sci. Eng. 139(2001), pp. 66-77.

\bibitem{williams} M. M. R. Williams,
{\it The effect of random geometry on the criticality of a
multiplying system IV: transport theory}, Nucl. Sci. Eng.
143(2003), pp. 1-18.




\end{thebibliography}
\end{document}